\catcode`\@=11

\font\tenmsb=msbm10
\font\sevenmsb=msbm7 \font\fivemsb=msbm5  \newfam\msbfam
\textfont\msbfam=\tenmsb\scriptfont\msbfam=\sevenmsb
\scriptscriptfont\msbfam=\fivemsb

\def\hexnumber@#1{\ifnum#1<10 \number#1\else \ifnum#1=10 A\else\ifnum#1=11
 B\else\ifnum#1=12 C\else \ifnum#1=13 D\else\ifnum#1=14 E\else\ifnum#1=15
 F\fi\fi\fi\fi\fi\fi\fi}
 \def\msb@{\hexnumber@\msbfam}

\mathchardef\hbar="0\msb@7E

\def\Bbb{\ifmmode\let\next\Bbb@\else
\def\next{\errmessage{Use \string\Bbb\space only in math mode}}\fi\next}
\def\Bbb@#1{{\Bbb@@{#1}}} \def\Bbb@@#1{\fam\msbfam#1}

\catcode`\@=\active

\def\CR{\hbox{{$\cal R$}}}
\def\CM{\hbox{{$\cal M$}}}
\def\CZ{\hbox{{$\cal Z$}}}
\def\CC{\hbox{{$\cal C$}}}

\def\cg{\hbox{{\sl g}}} 

\def\C{{\Bbb C}}

\def\rcross{{\triangleright\!\!\!<}}
\def\tens{\mathop{\otimes}}
\def\dcross{{\bowtie}}
\def\<{{\langle}}
\def\>{{\rangle}}
\def\la{{\triangleright}}
\def\ra{{\triangleleft}}

\def\eps{{\epsilon}}
\def\Nat{{\rm Nat}}
\def\Ad{{\rm Ad}}

\def\id{{\rm id}}

\def\und#1{{\underline {#1}}}

\def\proof{\goodbreak\noindent{\bf Proof\quad}}
\def\endproof{{\ $\hbox{$\sqcup$}\llap{\hbox{$\sqcap$}}$}\bigskip }

\def\o{{}_{\scriptscriptstyle(1)}}
\def\t{{}_{\scriptscriptstyle(2)}}
\def\th{{}_{\scriptscriptstyle(3)}}
\def\fo{{}_{\scriptscriptstyle(4)}}
\def\fiv{{}_{\scriptscriptstyle(5)}}
\def\bo{{}^{\bar{\scriptscriptstyle(1)}}}
\def\bt{{}^{\bar{\scriptscriptstyle(2)}}}
\def\uo{{{}^{\scriptscriptstyle(1)}}}
\def\ut{{{}^{\scriptscriptstyle(2)}}}
\def\uth{{{}^{\scriptscriptstyle(3)}}}
\def\umo{{{}^{\scriptscriptstyle -(1)}}}
\def\umt{{{}^{\scriptscriptstyle -(2)}}}
\def\umth{{{}^{\scriptscriptstyle -(3)}}}

\def\nquad{{\!\!\!\!\!\!}}

\def\cmath#1{\[\begin{array}{c} #1 \end{array}\]}
\def\eqn#1#2{\begin{equation}#2\label{#1}\end{equation}}
\def\ceqn#1#2{\begin{equation}\label{#1}
\begin{array}{c}#2\end{array}\end{equation}}

\documentstyle[11pt]{article}
\newtheorem{lemma}{Lemma}[section]
\newtheorem{propos}[lemma]{Proposition}

\newtheorem{corol}[lemma]{Corollary}

\begin{document}\baselineskip 22pt

{\ }\qquad\qquad \hskip 4.3in Damtp/95-69
\vspace{.2in}

\begin{center} {\LARGE  QUANTUM DOUBLE FOR QUASI-HOPF ALGEBRAS}
\\ \baselineskip 13pt{\ }
{\ }\\
 S. Majid\footnote{Royal Society University Research Fellow and
Fellow of Pembroke College, Cambridge}\\ {\ }\\ Department of
Mathematics, Harvard University\\ Science Center, Cambridge MA
02138, USA\footnote{Work completed during my visit 1995+1996}\\ +\\
Department of Applied Mathematics \& Theoretical Physics\\
University of Cambridge, Cambridge CB3 9EW\\
\end{center}
\begin{center}
Final version to appear Lett. Math. Phys. (accepted January 1998).
\end{center}

\vspace{10pt}
\begin{quote}\baselineskip 13pt
\noindent{\bf Abstract}
We introduce a quantum double quasitriangular quasi-Hopf algebra
$D(H)$ associated to any quasi-Hopf algebra $H$. The algebra
structure is a cocycle double cross product. We use categorical
reconstruction methods. As an example, we recover the quasi-Hopf
algebra of Dijkgraaf, Pasquier and Roche as the quantum double
$D^\phi(G)$ associated to a finite group $G$ and group $3$-cocycle
$\phi$.
\bigskip
\noindent Keywords: quantum double -- quasi-Hopf algebra -- finite
group -- cocycle -- category -- reconstruction.

\end{quote}
\baselineskip 19pt

\section{Introduction}

The quantum double\cite{Dri} of  V.G. Drinfeld is one of the most
important of
quantum group constructions. It associates to a Hopf algebra $H$ a
quasitriangular one. The quasitriangular structure leads to a
braiding in the
category of representations and many ensuing applications.

In this note we introduce the corresponding construction for
quasi-Hopf
algebras. Quasi-Hopf algebras have coproducts which are coassociative
only up
to a 3-cocycle $\phi\in H\tens H\tens H$\cite{Dri:qua}. This greater
freedom
allows, for example, the simplest formulation of quantum groups
$U_q(\cg)$ as
ordinary enveloping algebras $U(\cg)$ equipped with a quasi-Hopf
structure
obtained from the Knizhnik-Zamolodchikov equations\cite{Dri:qua}.
This is the
form in which quantum groups `naturally arise' in conformal field
theory and in
the theory of Vassilyev invariants, for example.

At first, it would appear difficult to define directly the
corresponding quantum double because, in Drinfeld's construction,
the quantum double $D(H)$ contains $H,H^{*\rm op}$ in a symmetrical
way; if $H$ is quasi-coassociative then $H^{*\rm op}$ is
quasi-associative and one might expect the double to be some kind
of hybrid object. We will see that this is not necessary. Instead,
we use a categorical formulation of the quantum double which must
necessarily give us as the double an ordinary quasi-Hopf algebra.

The categorical picture of the usual quantum double was also
provided by V.G. Drinfeld, in terms of the `center' or `double'
$\CZ(\CC)$ of a monoidal category. The same (slightly more general)
construction was introduced at the same time in \cite{Ma:rep} as a
generalised `Pontryagin dual' $\CC^\circ$ of a monoidal category.
Drinfeld observed\cite{Dri:pri} that this category is braided and
is ${}_{D(H)}\CM$  when $\CC={}_H\CM$, the modules over a Hopf
algebra. We recall the basic facts in the following Preliminaries
section. In Section~2, we compute the double category when $H$ is
the category of modules over a quasi-Hopf algebra $H,\phi$. We can
then use the Tannaka-Krein reconstruction in the generalised form
cf\cite{Ma:tan} to define $D(H)$ for quasi-Hopf algebras, in
Section~3. We cover in particular the example when $H=k(G)$ the
algebra of functions on a finite group $G$, equipped with a group
cocycle $\phi\in Z^3(G)$. The double in this case recovers the
quasi-Hopf algebra introduced by direct means in \cite{DPR:qua} in
connection with a `toy model' of conformal field
theory\cite{DijWit:top}.

\subsection*{Preliminaries}

We work over a general field $k$ or, with suitable care, over a
commutative
ring. Following Drinfeld\cite{Dri:qua}, a quasi-Hopf algebra means
$(H,\Delta,\eps,\phi,S,\alpha,\beta)$ where $H$ is a unital algebra,
$\Delta:H\to H\tens H$ and $\eps:H\to k$ are algebra homomorphisms
obeying
\eqn{qcoassoc}{(\id\tens\Delta)\circ\Delta =\phi
((\Delta\tens\id)\circ\Delta(\
))\phi ^{-1},\quad
(\id\tens\eps)\circ\Delta=\id=(\id\tens\eps)\circ\Delta,}
where $\phi\in H^{\tens 3}$ is invertible and a 3-cocycle in the
sense
\eqn{3-cocycle}{{}\nquad(1\tens \phi)((\id\tens\Delta\tens
\id)\phi)(\phi\tens
1)=((\id\tens \id\tens \Delta)\phi)((\Delta\tens \id\tens \id)\phi)}
and $(\id\tens\eps\tens\id)\phi=1\tens 1$. This part defines a
quasi-bialgebra.
In addition, we require $S:H\to H$ and $\alpha,\beta\in H$ such that
\ceqn{qant}{ \sum (Sh\o)\alpha h\t=\eps(h)\alpha,\quad \sum h\o\beta
S
h\t=\eps(h)\beta,\qquad\forall h\in H\\
\sum \phi \uo\beta(S\phi \ut)\alpha\phi \uth =1,\qquad
\sum (S\phi \umo )\alpha\phi \umt  \beta S\phi ^{-(3)}=1.}
We use the notation $\Delta h=\sum h\o\tens h\t$
and $\phi =\sum\phi \uo\tens\phi \ut\tens\phi \uth $. Similarly for
$\phi^{-1}$.

A quasi-bialgebra or quasi-Hopf algebra $H$ is quasitriangular if
there is an
invertible element $\CR\in H\tens H$ such that
\eqn{qqua}{ \nquad\ (\Delta\tens\id)\CR=\phi _{312}\CR_{13}\phi
_{132}^{-1}\CR_{23}\phi,\ \,  (\id\tens\Delta)\CR=\phi
_{231}^{-1}\CR_{13}\phi
_{213}\CR_{12}\phi ^{-1}}
in another standard notation. Explicitly,  $\phi _{213}=\sum \phi
\ut\tens\phi
\uo\tens\phi \uth $, etc.

A monoidal category means a category $\CC$ with objects $V,W,Z$ etc.,
a functor
$\tens:\CC\times\CC\to \CC$ equipped with an associativity natural
transformation consisting of functorial isomorphisms
$\Phi_{V,W,Z}:(V\tens
W)\tens Z\to V\tens (W\tens Z)$ obeying a pentagon
identity\cite{Mac:cat}.
There is also a compatible unit object and associated functorial
isomorphisms.
A braided category is a monoidal category equipped with a
commutativity natural
transformation consisting of functorial isomorphisms
$\Psi_{V,W}:V\tens W\to
W\tens V$ compatible with the unit and associativity structures in a
natural
way \cite{JoyStr:bra}. We generally suppress writing the tensor
product by
identity morphisms, as well as  isomorphisms associated with the
identity
object.

A {\em representation}\cite{Ma:rep} of a monoidal category $\CC$ in
itself is
an object $V$ of $\CC$ and a natural equivalence
$\lambda_V\in\Nat(V\tens\id,\id\tens V)$, i.e. functorial
isomorphisms
$\lambda_{V,W}:V\tens W\to W\tens V$ (functorial in $W$), such that
\eqn{repcat}{ \lambda_{V,\und 1}=\id,\qquad
\lambda_{V,Z}\circ\Phi_{W,V,Z}\circ\lambda_{V,W}=\Phi_{W,Z,V}\circ
\lambda_{V,W\tens Z}\circ\Phi_{V,W,Z},\qquad\forall W,Z\in\CC.}
The `Pontryagin dual' monoidal category $\CC^\circ$, or `double', has
as
objects such representations. Morphisms $\phi:(V,\lambda_V)\to
(W,\lambda_W)$
are morphisms $\phi:V\to W$ of $\CC$  obeying
\ceqn{morrep}{ (\id\tens
\phi)\circ\lambda_{V,Z}=\lambda_{W,Z}\circ(\phi\tens\id).}
The monoidal structure $(V,\lambda_V)\tens(W,\lambda_W)$ consists of
$V\tens W$
and the natural transformation
\ceqn{tensrep}{   \lambda_{V\tens W,Z}=\Phi_{Z,V,W}\circ\lambda_{V,Z}
\circ\Phi^{-1}_{V,Z,W}
\circ\lambda_{W,Z}\circ\Phi_{V,W,Z},\quad\forall Z\in\CC.}
The associator $\Phi_{V,W,Z}$ is the underlying one for the
category $\CC$, viewed as functorial isomorphisms in $\CC^\circ$.
The construction also works more generally for representations in
another monoidal category. In the present case, there is a braiding
\ceqn{brarep}{\Psi_{(V,\lambda_V),(W,\lambda_W)}=\lambda_{V,W}}
due to \cite{Dri:pri}. Other notations for this category are $D(\CC)$
or
$\CZ(\CC)$. Further results are in \cite{Ma:cat}.

If $H$ is a quasi-bialgebra or quasi-Hopf algebra, we denote by
${}_H\CM$ its
category of modules.
This forms a monoidal category with tensor product defined via
$\Delta$ and
with associativity transformation
\eqn{repqHa}{\Phi_{V,W,Z}((v\tens w)\tens z)=\sum \phi\uo\la v\tens
(\phi\ut\la w\tens \phi\uth\la z),\quad \forall v\in V,\ w\in W,\
z\in Z} and (in the quasitriangular case) braiding defined by
\eqn{repqHb}{\Psi_{V,W}(v\tens w)=\sum \CR\ut\la
w\tens
\CR\uo\la v,\quad \forall v\in V, w\in W.}
The forgetful functor is multiplicative\cite{Ma:tan} but not
monoidal unless $H$ is twisting equivalent to an ordinary Hopf
algebra. Further details are in \cite{Ma:book}.

\section{Double category of modules over a quasi-Hopf algebra}

Let $H,\phi$ be a quasi-Hopf algebra with bijective antipode and
$\CC={}_H\CM$.
We denote the action of $H$ on a module $V$ by $\la$.

\begin{lemma} Let $V\in\CC$. Natural transformations $\Nat(V\tens
\id,\id\tens
V)$ are in 1-1 correspondence with linear maps $V\to H\tens V$,
denoted
$v\mapsto  \sum v\bo\tens v\bt$, such that
\[\sum h\o v\bo\tens h\t\la v\bt=\sum (h\o\la v)\bo h\t\tens (h\o\la
v)\bt.\]
\end{lemma}
\proof Let $H_L$ denoted $H\in {}_H\CM$ by the left regular
representation. If
$\lambda_V$ is a natural transformation, the corresponding linear map
$\beta:V\to H\tens V$ is
\[\beta(v)=\lambda_{V,H_L}(v\tens 1).\]
For and $W\in \CC$ and $w\in W$ we consider the morphism $i_w:H_L\to
W$ defined
by $h\mapsto h\la w$. Then naturality of $\lambda_V$ implies that
$\lambda_{V,W}(v\tens w)=\lambda_{V,W}(v\tens i_w.1)=v\bo\la w\tens
v\bt$.
In particular, for any $h\in H$ we consider right-multiplication
$R_h:H_L\to
H_L$ as a morphism and hence $\lambda_{V,H_L}(v\tens
h)=\lambda_{V,H_L}(v\tens
R_h.1)=v\bo h\tens v\bt$. Then the assumption that each
$\lambda_{V,H_L}$ are
morphisms implies that
\[(h\o\la v)\bo h\t\tens (h\o\la v)\bt=\lambda_{V,H_L}(h\o\la v\tens
h\t.1)=h\o
v\bo\tens h\t \la v\bt,\]
which is the condition stated. Conversely, given $\beta$ obeying this
condition, we define $\lambda_{V,W}(v\tens w)=v\bo\la w\tens v\bt$
and verify
easily that this is a natural transformation.  \endproof

\begin{propos} Representations of ${}_H\CM$ in itself are in 1-1
correspondence
with  pairs $(V,\beta_V)$ where $V$ is an $H$-module and
$\beta_V:V\to H\tens
V$ obeys the condition in Lemma~2.1 and
\[ \phi\uo v\bo\tens (\phi\ut\la v\bt)\bo \phi\uth\tens(\phi\ut\la
v\bt)\bt=\phi\la\left((\phi\uo\la v)\bo\o \phi\ut\tens (\phi\uo\la
v)\bo\t
\phi\uth\tens (\phi\uo\la v)\bt\right),\]
where $\phi$  acts on $H_L\tens H_L\tens V$. We also require
$(\eps\tens\id)\circ\beta_V=\id$. The category $({}_H\CM)^\circ$
consists of
such objects and morphisms which intertwine the $H$ action and the
corresponding $\beta$. The category has  monoidal product
$(V,\beta_V)\tens
(W,\beta_W)$ built on $V\tens W$ as an $H$-module and
\[ \beta_{V\tens W}=\phi\la\left((\phi\umo\phi\uo\la v)\bo \phi\umt
(\phi\ut\la w)\bo \phi\uth\tens (\phi\umo\phi\uo\la v)\bt\tens
\phi\umth\la (\phi\ut\la w)\bt\right),\] where $\phi$ acts on
$H_L\tens V\tens W$. The associativity isomorphisms are given by
the same formula (\ref{repqHa}) as for ${}_H\CM$. The category is
braided, with
\[ \Psi_{V,W}(v\tens w)=v\bo\la w\tens v\bt.\]
\end{propos}
\proof We write out (\ref{repcat}) and (\ref{tensrep}) using the
identifications in Lemma~2.1. In the converse direction, the
$\lambda_{V,W}$
defined by $\beta_V$ are isomorphisms since $H$ using the assumed
inverse
antipode (this corresponds to existence of left duals by
\cite{Ma:rep}).
\endproof

These steps are similar to the computation of the double category
for ordinary Hopf algebras, with the additional presence of $\phi$
in our quasi-Hopf algebra case. Whereas one could also come to such
a category of `crossed modules'  as a generalisation of Whitehead's
crossed $G$-sets\cite{Whi:com}, one really needs the above double
approach in the quasi-Hopf algebra case, in order to place the
$\phi$ correctly.

\section{Double of a quasiHopf algebra}

We are now in a position to define the double quasi-Hopf algebra
$D(H)$ as that obtained by Tannaka-Krein reconstruction from the
category constructed in the preceding section as the automorphisms
of the forgetful functor.

\begin{corol} If $(H,\phi)$ is a finite-dimensional quasi-Hopf algebra, it
has a quantum double $D(H)$ uniquely defined, up to isomorphism, as
a quasitriangular quasiHopf algebra such that its category of
representations is $({}_H\CM)^\circ$ in Proposition~2.2. In
particular, it may be built on the vector space $H^*\tens H$ with
$H$ as a sub-quasiHopf algebera.
\end{corol}
\proof The forgetful functor from $({}_H\CM)^\circ$ to $Vec$ is
multiplicative
and hence we can use the reconstruction theorem \cite{Ma:tan}. This
builds $D(H)$ such that we can identify
${}_{D(H)}\CM=({}_H\CM)^\circ$ as braided monoidal categories. It
is clear from the characterisation of the latter in Proposition~2.2
that $D(H)$ may be built on $H^*\tens H$ as a vector space with a
certain product and coproduct: the action of $(f\tens h)\in D(H)$
is
\eqn{douact}{ (f\tens h)\la v=\<f,(h\la v)\bo \> (h\la v)\bt}
where $(V,\la,\beta)$ is the corresponding object of
$({}_H\CM)^\circ$ and $\<\ ,\ \>$ is the evaluation pairing.
Moreover, the forgetful functor factors through the monoidal
functor $({}_H\CM)^\circ\to {}_H\CM$, which corresponds to an
inclusion $H\subset D(H)$ as $h\equiv 1\tens h$. Note that the form
of $\Phi$ gives immediately
\eqn{douphi}{ \phi_{D(H)}=\phi}
under this inclusion. Likewise, the form of $\Psi$ in
Proposition~2.2 and comparison with (\ref{repqHb}) gives
immediately
\eqn{douR}{ \CR=\sum_a (f^a\tens 1)\tens (1\tens e_a),}
where $\{e_a\}$ is a basis of $H$ and $\{f^a\}$ is a dual basis.
The reconstructed product and coproduct are more complex, although
reducing to those for the usual quantum double $H^{*\rm op}\dcross
H$ when $\phi=1$.
\endproof

If one wants explicit formulae for the product and coproduct of
$D(H)$, they are immediately obtained from the formulae in
Lemma~2.1 and Proposition~2.2. We can write the former as
\[ (\phi\umo h\la v)\bo \phi\umt\beta S\phi\umth\tens
(\phi\umo h\la v)\bt
=h\o\o (\phi\umo\la v)\bo \phi\umt\beta S(\phi\umth h\t)
\tens h\o\t\la (\phi\umo\la v)\bt\]
in view of the quasi-coassociativity and antipode axioms for $H$.
Applying $\<f,\ \>$ for $f\in H^*$ to the first factor, and making
use of (\ref{douact}), we eliminate $v$ and obtain the relations
\eqn{doufh}{ f\o\cdot \phi\umo\cdot h\<f\t,\phi\umt\beta S\phi\umth\>
=\<f\o,h\o\o\>\<f\th,\phi\umt\beta S(\phi\umth h\t)\> h\o\t\cdot
f\t\cdot \phi\umo} for all $h\in H$ and $f\in H^*$. Here we write
$h\o\o\equiv 1\tens h\o\o$ and $f\o\equiv f\o\tens 1$, etc., and
$\cdot$ is the product in $D(H)$. Moreover, $f\o\tens f\t$ etc.,
denotes the coassociative coproduct of $H^*$ dual to the product of
$H$. Similarly, applying $\<g,\
\>\tens\<f,\ \>$ to the quasi-coaction condition for
$\beta_V$ in Proposition~2.2,
using (\ref{douact})
and cancelling $v$ gives immediately
\eqn{doufg}{ \<g\o,\phi\uo\>\<f\t,\phi\uth\> f\o\cdot \phi\ut\cdot g\t
=\<g\o,\phi'{}\uo\>\<f\o,\phi'{}\ut\>\<g\th,\phi\ut\> \<f\th,\phi\uth\>
\phi'{}\uth\cdot(g\t f\t)
\cdot\phi\uo}
for all $f,g\in H^*$. Here $\phi'$ denotes a second copy of $\phi$
and $g\t f\t$ is multiplied in the (non-associative) product on
$H^*$ dual to the coproduct of $H$. Finally, applying $\<f,\
\>$ to the formula for $\beta_{V\tens W}$ in Proposition~2.2 and
cancelling $v,w$ gives
\eqn{doudelta}{ \Delta_{D(H)}f=\<f\o,\phi'\uo\>\<f\th,\phi\umt\>
\<f\fiv,\phi\uth\>
\phi'\ut\cdot f\t\cdot \phi\umo\phi\uo\tens \phi'\uth\phi
\umth\cdot f\fo\cdot\phi\ut.}
for $f\in H^*$. We already know that $\Delta_{D(H)}h=\Delta h$ for
$h\in H$. These more explicit formulae
(\ref{douphi})--(\ref{doudelta}) for $D(H)$ correspond directly to
the characterisation of its representations in Lemma~2.1 and
Proposition~2.2.\footnote{We have added them here at the request of
the referee. We note that in the meantime the recent preprint
`Doubles of quasi-quantum groups' by F. Hausser and F. Nill,
following up the preprint version of the present paper, provides
some further explicit formulae for the above $D(H)$.}

To illustrate this theory, we content ourselves with the simplest
case, which is, however, the case relevant for conformal field
theory so far. Thus, let $G$ be a finite group and $\phi\in
Z^3(G)$, i.e. a 3-cocycle in the sense
\[ \phi(y,s,t)\phi(x,ys,t)\phi(x,y,s)=\phi(x,y,st)\phi(xy,s,t),\quad
\phi(x,e,y)=1\]
for all $x,y,s,t\in G$ and $e$ the group identity element. Let $k(G)$
be the
Hopf algebra of functions on $G$ with coproduct $(\Delta
f)(x,y)=f(xy)$. We
view it as a quasi-Hopf algebra $k^\phi(G)$ with $\phi\in k(G)^{\tens
3}$. We
also make use of $kG$, the group algebra of $G$.

\begin{propos} An object of the double category
$({}_{k^\phi(G)}\CM)^\circ$  is
 a $G$-graded vector space $V$ and a right cocycle action $\ra$ of
$kG$ on $V$
which is compatible with the grading, in the sense
\[ (v\ra x)\ra y
={\phi(x,y,y^{-1}x^{-1}|v|xy)\phi(|v|,x,y)\over\phi(x,x^{-1}|v|x,y)}
v\ra(xy),\quad v\ra 1=v,\quad |v\ra x|=x^{-1}|v|x\]
for all $v\in V$ homogeneous of degree $|v|$ and $x,y\in G$.
Morphisms in the
category are linear maps preserving the right action and the grading.
The
tensor product grading  and  action are given by
\[ |v\tens w|=|v||w|,\quad (v\tens w)\ra
x={\phi(x,x^{-1}|v|x,x^{-1}|w|x)\phi(|v|,|w|,x)\over
\phi(|v|,x,x^{-1}|w|x)}
v\ra x\tens w\ra x,\]
for homogeneous $v,w$. We call this the category of  cocycle crossed
$G$-modules. It is braided, with $\Psi_{V\tens W}(v\tens w)=w\ra
|v|\tens v$.
\end{propos}
\proof  If $V$ is a $k(G)$-module, the projection operators given by
the left
action of Kronecker delta-functions $\delta_x$ provide a
decomposition of $V$
into subspaces of degree $x$. Conversely, if $V$ is $G$-graded, it
becomes a
$k(G)$-module by $f\la v=f(|v|)v$ on homogeneous elements. Next, we
write
$\beta_V$ in Proposition~2.2 as a right `quasi-action' by $kG$
according to the
correspondence $v\ra x=\sum v\bo(x) v\bt$. The condition in
Proposition~2.2 in
these terms comes out as stated. It is easy to see that $\chi:G\times
G\to
k(G)$ defined by
\eqn{chi}{
\chi(x,y)(s)={\phi(x,y,y^{-1}x^{-1}sxy)\phi(s,x,y)\over\phi(x,x^{-1}s
x,y)} } is a 2-cocycle in $Z^2_{\Ad}(G,k(G))$ with values in $k(G)$
as a left module induced by the group adjoint action. This is the
sense in which $\ra$ is a cocycle action. The remaining structure
easily computes in this case as stated.
\endproof

This category generalises J.H.C. Whitehead's notion of crossed
$G$-sets\cite{Whi:com} to the case of a non-trivial 3-cocycle $\phi$.
It is
easy to identify the algebra, $D^\phi(G)$, say, with the
representations of
which this category can be identified. Namely, cocycle
representations are
naturally identified with modules of the corresponding cocycle cross
product
algebra.

\begin{propos} The quasi-Hopf algebra double in Corollary~3.1 reduces in
this
example to the quasi-Hopf algebra $D^\phi(G)=kG^{\rm op}\rcross_\chi
k(G)$ in
\cite{DPR:qua}. Explicitly,
\cmath{ (x\tens\delta_s)\cdot(y\tens\delta_t)=yx\tens\delta_t
\delta_{ysy^{-1},t}\chi(y,x)(t)\\
\Delta(x\tens\delta_s)=\sum_{ab=s}{\phi(x,x^{-1}ax,x^{-1}bx)\phi(a,b,
x)\over
\phi(a,x,x^{-1}bx)}x\tens\delta_a\tens x\tens\delta_b,\quad
\eps(x\tens\delta_s)=\delta_{s,e}}
and $\phi\in D^\phi(G)^{\tens 3}$ by the standard inclusion. There is
a
quasitriangular structure
$\CR=\sum (\delta_x\tens 1)\tens (1\tens x)$.
\end{propos}
\proof This is a special case of Section~3. However, it is also easy
enough to
see directly in the present example. We view
$\chi_{21}(x,y)=\chi(y,x)$ as a
right-handed 2-cocycle on $G^{\rm op}$ with values in $k(G)$ viewed
as  a right
$G^{\rm op}$-module (by the adjoint action). The 2-cocycle property
means
$\chi_{21}(x,y)\ra s\chi_{21}(xy,s)=\chi_{21}(y,s)\chi_{21}(x,ys)$.
Hence by a
standard construction for cocycle cross product algebras by group
cocycle-actions\cite{Swe:coh}, we have an algebra $kG^{\rm
op}\rcross_\chi
k(G)$ with product $(x\tens f)(y\tens h)=yx\tens \chi(y,x)(f\ra y)h$
for
$f,h\in k(G)$. This has the form shown on delta-functions. Moreover,
the
category in Proposition~4.1 can be identified with its left modules
in the
obvious way. If $V$ is an object of the category, the corresponding
left module
is $(x\tens f)\la v=v\ra x f(|v|)$ on homogeneous elements. Given
this
identification, the tensor product of objects in Proposition~4.1
corresponds to
the map coproduct shown.
\endproof

This is the content of our more formal Tannaka-Krein arguments in
the finite groups setting. It is clear that the result is a special
case of the general formulae (\ref{douphi})--(\ref{doudelta}) of
$D(H)$. It should also be clear that this double quasi-Hopf algebra
construction has an immediate generalisations to the case of an
infinite-dimensional quasi-Hopf algebra and a dual quasi-Hopf
algebra paired with it. (And we do not really require  an antipode,
provided the pairing is convolution-invertible in the usual sense.)
Alternatively, one may work over $\C[[\hbar]]$ and take suitable
topological duals as in the usual quantum double construction. In
such a setting, an interesting project for further work would be to
apply the double quasi-Hopf algebra construction to
$H=(U\cg,\phi)$, where $\cg$ is a complex semisimple Lie algebra
and $\phi$ is Drinfeld's 3-cocycle obtained by solving the KZ
equation\cite{Dri:qua}. For example, when $g=su_2$, this provides
in principle a quasi-Hopf algebra approach to the q-Lorentz group.


\begin{thebibliography}{10}

\bibitem{Dri}
V.G. Drinfeld.
\newblock Quantum groups.
\newblock In A.~Gleason, editor, {\em Proceedings of the {ICM}},
pages
  798--820, Rhode Island, 1987. AMS.

\bibitem{Dri:qua}
V.G. Drinfeld.
\newblock Quasi{H}opf algebras.
\newblock {\em Leningrad Math. J.}, 1:1419--1457, 1990.

\bibitem{Ma:rep}
S.~Majid.
\newblock Representations, duals and quantum doubles of monoidal
categories.
\newblock {\em Suppl. Rend. Circ. Mat. Palermo, Ser. II},
26:197--206, 1991.

\bibitem{Dri:pri}
V.G. Drinfeld.
\newblock Private communication, {F}ebruary, 1990.

\bibitem{Ma:tan}
S.~Majid.
\newblock Tannaka-{K}rein theorem for quasi{H}opf algebras and other
results.
\newblock {\em Contemp. Maths}, 134:219--232, 1992.

\bibitem{DPR:qua}
V.~Pasquier R.~Dijkgraaf and P.~Roche.
\newblock Quasi-quantum groups related to orbifold models.
\newblock In {\em Proc. Modern Quantum Field Theory, Tata Institute,
Bombay,
  1990}, pages 375--383.

\bibitem{DijWit:top}
R.~Dijkgraaf and E.~Witten.
\newblock Topological gauge theories and group cohomology.
\newblock {\em Commun. Math. Phys}, 129:393--429, 1990.

\bibitem{Mac:cat}
S.~Mac Lane.
\newblock {\em Categories for the Working Mathematician}.
\newblock Springer, 1974.
\newblock GTM vol. 5.

\bibitem{JoyStr:bra}
A.~Joyal and R.~Street.
\newblock Braided monoidal categories.
\newblock Mathematics Reports 86008, Macquarie University, 1986.

\bibitem{Ma:cat}
S.~Majid.
\newblock Braided groups and duals of monoidal categories.
\newblock {\em Canad. Math. Soc. Conf. Proc.}, 13:329--343, 1992.

\bibitem{Ma:book}
S.~Majid.
\newblock {\em Foundations of Quantum Group Theory}.
\newblock Cambridge University Press, 1995.

\bibitem{Whi:com}
J.H.C. Whitehead.
\newblock Combinatorial homotopy, {II}.
\newblock {\em Bull. Amer. Math. Soc.}, 55:453--496, 1949.

\bibitem{Swe:coh}
M.E. Sweedler.
\newblock Cohomology of algebras over {H}opf algebras.
\newblock {\em Ann. Math.}, pages 205--239, 1968.

\end{thebibliography}

\end{document}